# LARGE DEVIATIONS PROBLEMS FOR STAR NETWORKS: THE MIN POLICY[1]

By Franck Delcoigne and Arnaud de La Fortelle[2]

*EDF R&D and INRIA*

We are interested in analyzing the effect of bandwidth sharing for telecommunication networks. More precisely, we want to calculate which routes are bottlenecks by means of large deviations techniques. The method is illustrated in this paper on a star network, where the bandwidth is shared between customers according to the so-called min policy. We prove a sample path large deviation principle for a rescaled process $n^{-1}Q_{nt}$, where $Q_t$ represents the joint number of connections at time $t$. The main result is to compute the rate function *explicitly*. The major step consists in deriving large deviation bounds for an *empirical generator* constructed from the join number of customers and arrivals on each route. The rest of the analysis relies on a suitable change of measure together with a localization procedure. An example shows how this can be used practically.

**1. Introduction.**

*The model.* Consider a star shaped network (see Figure 1) consisting of $N$ links connected to the other $N-1$ through a central hub: there are $N(N-1)/2$ routes of length two. In the sequel the set of links (or channels) is denoted by $\mathcal{S} = \{1, \ldots, N\}$, whereas the set of routes is simply the set of unordered
two-uples $ij$, $i, j \in \mathcal{S}$: for the sake of simplicity, we do not distinguish between $ij$ and $ji$ (i.e., we consider nonoriented routes), but there is no additional difficulty to handle oriented routes. Denote by $q_{ij}(t)$ [resp. $q_i(t)$] the number of calls (or connections) on route $ij$ (resp. the number of calls involving channel $i$) at time $t$. Each link has a capacity (or bandwidth) equal to $C_i$ (expressed,

Received February 2002; revised May 2003.
[1]This work is a part of contract 981B016 between INRIA and France Telecom R&D.
[2]Supported by the French Ministère de l'Équipement.
*AMS 2000 subject classifications.* Primary 60F10; secondary 60K30.
*Key words and phrases.* Large deviations, rate function, empirical generator, change of measure, contraction principle, entropy, star network, bandwidth sharing, min protocol.

e.g., in bits per second in the context of communication networks). Note that $q_i(t) = \sum_j q_{ij}(t)$. Then $Q(t,x) = (q_{ij}(t), i,j \in \mathcal{S})$ represents the state of the network at time $t$ when it starts initially from state $x$. For the sake of simplicity, we shall sometimes omit $x$ or $t$ when they do not play a role.

Documents to be transferred arrive on route $ij$ according to a Poisson process of rate $\lambda_{ij}$. We shall denote by $\mathcal{R}$ the set of active routes, that is, with $\lambda_{ij} > 0$. The size of a document (expressed in bits) on route $ij$ is supposed to be exponentially distributed with parameter $\mu_{ij}$. Each document on route $ij$ is allocated a portion $\nu_{ij}(x)/x_{ij}$ of the bandwidth when the state of the network is $x$. Hence, a document on route $ij$ is transferred at rate $\mu_{ij}\nu_{ij}(x)$. There are several possibilities in order to allocate a fair proportion of the bandwidth to customers. A classical one is to choose the coefficients $\nu_{ij}(x)$ according to the max–min fairness allocation.

The star network is proposed as a model for a router where the bandwidth is shared fairly between the different connections. However, the max–min fairness allocation is not explicit and hard to analyze at first. In order to get a more tractable model, we focus on the min policy,

$$\nu_{ij}(x) = \begin{cases} x_{ij}\dfrac{C_i}{x_i} \wedge \dfrac{C_j}{x_j}, & \text{if } x_{ij} > 0, \\ 0, & \text{otherwise.} \end{cases}$$

It has been shown in Fayolle, de La Fortelle, Lasgouttes, Massoulie and Roberts (2001) that the system under the max–min fairness allocation is stochastically smaller than the one with the min policy and that the network is ergodic if, and only if,

$$(1.1) \qquad \sum_j \frac{\lambda_{ij}}{\mu_{ij}} < C_i \qquad \forall i \in \mathcal{S}.$$

However, it appears very difficult to compute quantities of interest like the mean transfer time of a document, so we turn to asymptotic analysis, especially large deviations.

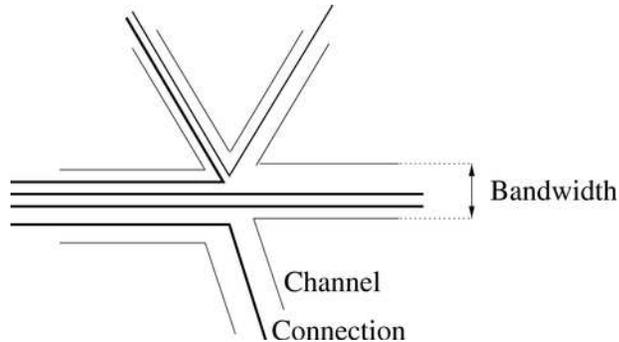

Fig. 1. *The (asymmetric) star network.*

*Previous work.* Lots of work has been devoted to the analysis of telecommunication networks. The model (star network and min policy) is described within telecommunication context in Fayolle, de La Fortelle, Lasgouttes, Massoulie and Roberts (2001). In the present paper we aim at deriving a sample path LDP for the rescaled process

$$Q_x^n \stackrel{\text{def}}{=} \left\{ \frac{1}{n} Q(nt, [nx]), \ t \geq 0 \right\}.$$

Our main concern is to identify explicitely the rate function. This is a preliminary step in order to obtain large deviation bounds in stationary regime. This issue is discussed in Section 2.

The major difficulty comes from the fact that the coefficients of the generator are not spatially continuous [the service rate $\mu_{ij}(x)$]. It seems that one of the first papers dealing with large deviations for processes with discontinuous statistics is Dupuis, Ishii and Soner ((1990)), where the case of Jackson networks was investigated using partial differential equations techniques. In Dupuis and Ellis ((1995)) a sample path LDP is proved for a wide class of jump Markov processes with discontinuous statistics. However, the methodology of proof uses subadditivity arguments and the rate function is not identified; moreover, there is a uniform reachability condition that our model does not fit. The identification of the rate function in this general framework is still an open problem when the dimension of the network is arbitrary. General results were obtained in Dupuis and Ellis ((1992)) and Ignatyuk, Malyshev and Shcherbakov ((1994)), where the LDP has been established. Nevertheless, in such examples, there are at most two boundaries with codimension one or two where discontinuity arises. Using special features of the models and the fact that fluid limits could be completely identified, this program was carried out, for example, in Atar and Dupuis ((1999)), Ignatiouk-Robert ((2000)) and in Delcoigne and de La Fortelle ((2002)).

*Structure of the paper.* An example (Section 2) shows how the rate function expression can be used to compute decay rate for tails of stationary distribution. In Section 3 we introduce the central notion of *localized model* and of *empirical generator*; the rate functions (local and global) are studied. In Section 4 the local LDP is proved by mean of a change of measure and the identification of the local rate function is worked out for ergodic networks. In Section 5 the sample path LDP is stated. In Section 6 we get rid of the ergodicity assumption: in our opinion this is the main contribution of the present paper since the methodology used allows a treatment of more complex and realistic protocols like max–min-fair. This issue will be investigated in a forthcoming paper.

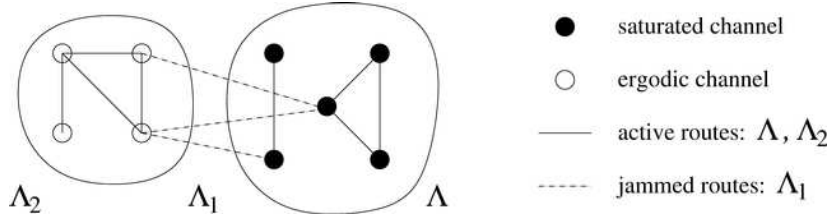

FIG. 2. *Representation of a star-shaped network: lines symbolize routes using two channels (circles at the ends of the lines). The routes are partitioned into saturated routes ($\Lambda$), jammed routes ($\Lambda_1$)—the service rate being null on these routes since all channels belonging to $\Lambda$ are saturated—and ergodic routes ($\Lambda_2$).*

*Notation.* In our settings, $\{Q(t, x_0),\ t \geq 0\}$ is a Markov process with generator $R$ such that, for all bounded real function $f$ on $\mathbb{Z}_+^{\mathcal{R}}$,

$$Rf(x) = \sum_{y \in \mathbb{Z}_+^{\mathcal{R}}} q(x,y)(f(y) - f(x)) \qquad \forall x \in \mathbb{Z}_+^{\mathcal{R}},$$

where

$$q(x,y) \stackrel{\text{def}}{=} \begin{cases} \lambda_{ij}, & \text{if } y - x = e_{ij}, \\ \mu_{ij}(x) \stackrel{\text{def}}{=} \mu_{ij} x_{ij} \dfrac{C_i}{x_i} \wedge \dfrac{C_j}{x_j}, & \text{if } y - x = -e_{ij}, \\ 0, & \text{otherwise,} \end{cases}$$

using the convention that $0/0 = 0$ (i.e., when $x_{ij} = 0$).

- For any set $A$, $A^c$ will denote its complementary and $\mathbb{1}_{\{A\}}$ its indicator function;
- $D([0,T], \mathbb{R}_+^{\mathcal{R}})$ is the space of right continuous functions $f : [0,T] \to \mathbb{R}_+^{\mathcal{R}}$ with left limits, endowed with the Skorokhod metric denoted by $d_d$;
- $\mathcal{C}([0,T], \mathbb{R}_+^{\mathcal{R}})$ is the space of continuous functions equipped with the metric of the uniform convergence denoted by $d_c$.

DEFINITION 1.1 (Face). For $x \in \mathbb{R}_+^{\mathcal{R}}$, the face $\Lambda(x)$ is defined by

$$\Lambda(x) \stackrel{\text{def}}{=} \{ij \in \mathcal{R} : x_{ij} > 0\}.$$

By an abuse of notation, we will also call face $\Lambda$

(1.2) $\qquad \{y \in \mathbb{R}_+^{\mathcal{R}} : y_{ij} > 0,\ \forall ij \in \Lambda,\ \text{and}\ y_{ij} = 0,\ \forall ij \in \Lambda^c\}.$

A partition of the routes (see Figure 2) is defined by $\Lambda$ and

$$\Lambda_1 \stackrel{\text{def}}{=} \{ij \in \Lambda^c : \exists k \in \mathcal{S},\ ik \in \Lambda \text{ or } jk \in \Lambda\},$$

$$\Lambda_2 \stackrel{\text{def}}{=} \{ij \in \Lambda^c : \forall k \in \mathcal{S},\ ik \notin \Lambda \text{ and } jk \notin \Lambda\}.$$

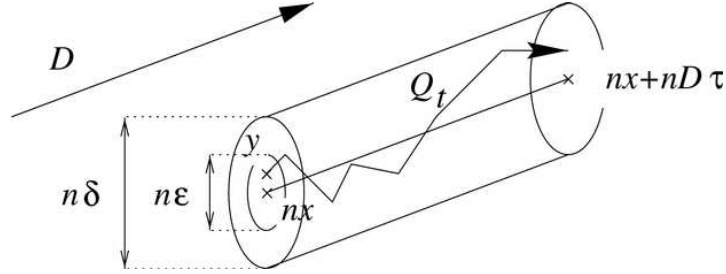

FIG. 3. *Structure of the local linear bounds of Theorem* 1.2. $L(x,D)$ *is the cost per unit time for the path $Q(t,y)$ (starting near $nx$) to stay in the neighborhood of $nx + Dt$ over the time $t \in [0, n\tau]$.*

The vector space relative to $\Lambda$ is defined by

$$\mathbb{R}^\Lambda \stackrel{\text{def}}{=} \{y \in \mathbb{R}^\mathcal{R} : y_{ij} = 0, \ \forall ij \in \Lambda^c\}.$$

*Results.* For ergodic networks, our main result is the local linear large deviation bounds of Theorem 1.2. The notation is explained by Figure 3.

THEOREM 1.2. *Assume that $Q$ is ergodic and let $x \in \mathbb{R}^\mathcal{R}_+$ and $D \in \mathbb{R}^{\Lambda(x)}$. Then, writing $\lim_{\tau,\delta,\epsilon \to 0}$ for $\lim_{\tau \to 0} \lim_{\delta \to 0} \lim_{\epsilon \to 0}$,*

$$\lim_{\tau,\delta,\epsilon \to 0} \inf_{|y-nx|<\epsilon n} \liminf_{n \to \infty} \frac{1}{n\tau} \log \mathbb{P}\left[\sup_{t \in [0,n\tau]} |Q(t,y) - nx - Dt| < \delta n\right]$$
(1.3)
$$= \lim_{\tau,\delta,\epsilon \to 0} \sup_{|y-nx|<\epsilon n} \limsup_{n \to \infty} \frac{1}{n\tau} \log \mathbb{P}\left[\sup_{t \in [0,n\tau]} |Q(t,y) - nx - Dt| < \delta n\right].$$

*Moreover, if $\Lambda$ and the drift $D \in \mathbb{R}^\Lambda$ are fixed, the preceding limit in $\tau$ is uniform w.r.t. to $x$ in compact sets of $\Lambda$ (see Definition 1.1). The common value of these limits is denoted by $-L(x,D)$ and*

$$(1.4) \qquad L(x,D) = \sum_{ij \in \Lambda(x) \cup \Lambda_1(x)} l(D_{ij} \| \lambda_{ij}, \mu_{ij}(x)),$$

*where*

$$(1.5) \quad l(D\|\lambda,\mu) \stackrel{\text{def}}{=} D \log\left(\frac{D + \sqrt{D^2 + 4\lambda\mu}}{2\lambda}\right) + \lambda + \mu - \sqrt{D^2 + 4\lambda\mu} \geq 0$$

*stands for the cost that a suitably normalized $M/M/1$ queue with parameters $\lambda$ and $\mu$, starting far from the origin, follows the drift $D$ [see, e.g., Shwartz and Weiss ((1995))].*

Let us explain briefly the meaning of the different terms appearing in $L(x, D)$ [see (1.4)]. Owing to the fact that the service rate $\mu_{ij}(x)$ tends to 0 when $x_{ij}$ becomes null, while $x_i \vee x_j$ remains strictly positive, the arrivals must be cut on the routes $ij \in \Lambda_1(x)$ in order to keep these routes in a neighborhood of 0. The cost to do this is $\sum_{ij \in \Lambda_1(x)} \lambda_{ij}$; indeed, $l(0\|\lambda_{ij}, 0) = \lambda_{ij}$. Since the arrivals are cut on the routes $ij \in \Lambda_1(x)$, the routes $ij \in \Lambda_2(x)$ are isolated from the rest of the network (see Figure 2) and so by (1.1) this set of routes behaves as an ergodic star network [with $\mathcal{R} = \Lambda_2(x)$] since $Q$ is ergodic by assumption. Hence, the cost for the components $ij \in \Lambda_2(x)$ to stay in a neighborhood of 0 is null. Now locally, the routes $ij \in \Lambda(x)$ behave as a set of independent $M/M/1$ queues with arrival and service rates $\lambda_{ij}$ and $\mu_{ij}(x)$. The corresponding terms in $L(x, D)$ represent the cost that this set of queues follows the prescribed drift $D$.

The proof is done introducing a functional so called empirical generator consisting of $Q_t$ and of the join number of arrivals on routes belonging to $\Lambda(x) \cup \Lambda_1(x)$. In Section 4 large deviation bounds are obtained for the localized empirical generator from which Theorem 1.2 is derived using an adaptation of the contraction principle.

Theorem 1.2 has been stated for ergodic networks. However, large deviations bounds can be obtained for transient networks, at the cost of some more detailed analysis. This is an important feature since it is linked with the study of networks under max–min-fair allocation (or similar ones). The reason is that, for an *ergodic* network under max–min-fair allocation, when some routes $ij \in \Lambda$ are made saturated (i.e., for localized models), the rest of the routes (in our notation $\Lambda^c$) can behave as a *transient* network, still under max–min-fair allocation: the local rate function must include the cost for this transient network to stay near 0. This is to the opposite of our framework, where only ergodic networks are considered, for which the cost to stay around 0 is null. However, our methodology allows one to compute cost for a transient network under the min policy to stay around 0 and as a corollary the rate function without ergodicity assumptions [see (6.6)]. The result is stated and discussed in Section 6.

Moreover, the topology of the network can be extended, as well as the length of the routes, (but not arbitrarily) to include more realistic networks. However, the notation becomes very heavy and our aim is to present tools [extending those developed for polling networks in Delcoigne and de La Fortelle (2002)] in a fairly simple way for achieving the above program.

Now, the rate function $I_T(\cdot)$ for the sample path LDP is expressed as

$$(1.6) \quad I_T(\varphi) \stackrel{\text{def}}{=} \begin{cases} \int_0^T L(\varphi(t), \dot\varphi(t))\, dt, & \text{if } \varphi \text{ is absolutely continuous,} \\ +\infty, & \text{otherwise.} \end{cases}$$

REMARK. $I_T(\cdot)$ is defined by all the values $L(x,D)$ with $x \in \mathbb{R}_+^{\mathcal{R}}$ and $D \in \mathbb{R}^{\Lambda}(x)$ (i.e., the values treated by Theorem 1.2). Indeed, assume that for some $t$, $\varphi_{ij}(t) = 0$ and $\dot{\varphi}_{ij}(t)$ exists. Since $\varphi_{ij}(t) \leq \varphi_{ij}(s)$ for all $s$, this implies $\dot{\varphi}_{ij}(t) \leq 0$. Then, necessarily, $\dot{\varphi}_{ij}(t) = 0$. Moreover, $\varphi$ being absolutely continuous, $\dot{\varphi}_{ij}(t)$ exists for almost all $t$.

Define the level set

$$(1.7) \qquad \Phi_x(K) \stackrel{\text{def}}{=} \{\varphi \in D([0,T], \mathbb{R}_+^{\mathcal{R}}) : I_T(\varphi) \leq K, \ \varphi(0) = x\}.$$

The final result is the following theorem.

THEOREM 1.3 (Sample path LDP). *Assume $Q$ is ergodic. The sequence $\{Q_x^n, \ n \geq 1\}$ satisfies an LDP in $D([0,T], \mathbb{R}_+^{\mathcal{R}})$ with good rate function $I_T(\cdot)$: for every $T > 0$, $x \in \mathbb{R}_+^{\mathcal{R}}$,*

(i) *for $C \subset \mathbb{R}_+^{\mathcal{R}}$ compact, $\bigcup_{x \in C} \Phi_x(K)$ is compact in $\mathcal{C}([0,T], \mathbb{R}_+^{\mathcal{R}})$;*
(ii) *for each closed set $F$ of $D([0,T], \mathbb{R}_+^{\mathcal{R}})$,*

$$\limsup_{n \to \infty} \frac{1}{n} \log \mathbb{P}[Q_x^n \in F] \leq -\inf\{I_T(\phi), \ \phi \in F, \phi(0) = x\};$$

(iii) *for each open set $O$ of $D([0,T], \mathbb{R}_+^{\mathcal{R}})$,*

$$\liminf_{n \to \infty} \frac{1}{n} \log \mathbb{P}[Q_{x,s}^n \in O] \geq -\inf\{I_T(\phi), \ \phi \in O, \phi(0) = x\}.$$

Its proof is discussed in Section 5.

**2. Example.** We would like to emphasize what kind of further results we aim at deducing from the sample path large deviation principle. First, it seems that the optimal paths of large deviation *cannot, in general,* be calculated, but some special solutions may be, leading to explicit expressions for the asymptotics of stationary distribution (which is not known). This is a performance criteria of practical value: bounds for buffer size could be optimized, or simulation accuracy (through importance sampling using the change of measures associated to optimal paths) could be improved.

Freidlin and Wentzell's works exposed in Freidlin and Wentzell ((1984)) suggest that the tail of the stationary distribution of the link $i$ is related to $I_T(\cdot)$ by the following formula:

$$(2.1) \quad \lim_{n \to \infty} \frac{1}{n} \log \mathbb{P}[q_i > n] = -\inf_{T \geq 0} \inf_{\varphi} \{I_T(\varphi) : \varphi(0) = 0, \varphi_i(T) = 1\}.$$

Although technical, it is reasonable to argue that the preceding equality can be checked in our case. However, this leads to an infinite-dimensional optimization problem. Nonetheless, by comparison with a processor sharing model, it seems that we can have partial information. Indeed, if the optimal path leading to the saturation of a channel $i$ is such that this channel is always the bottleneck (i.e., $C_i/x_i$ is minimal) then the service rate can be written, for each connection

$$\mu_{ij}(x) = \mu_{ij} C_i \frac{x_{ij}}{x_i}.$$

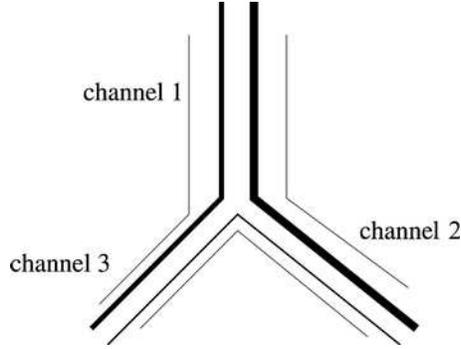

Fig. 4. *The network studied as example. The capacities of the respective channel are $C_1 = 3$, $C_2 = 2$ and $C_3 = 1$. The arrivals and document size are $\lambda_{12} = \mu_{12} = 1$, $\lambda_{23} = 1$, $\mu_{23} = 2$ and $\lambda_{13} = x$, $\mu_{13} = 1$.*

This is exactly the service rate for a processor sharing model which is a well-known model [e.g., the stationary distribution is explicit; see Baskett, Chandy, Muntz and Palacios (1975)]. With some calculations we can find a *necessary* condition for a channel $i$ to behave like a processor sharing (having thus the same stationary distribution decay rate). Otherwise there are more complex interactions between channels.

To illustrate the changes in the channels behavior with the load, we estimated the queues decay rates for different arrivals rates $\lambda_{13} \stackrel{\text{def}}{=} x$. We simulated the network described in Figure 4 and obtained statistics for the stationary queue length $\mathbb{P}[Q_i = n]$ decay rate. These results are compared with the decay rates of processor sharing models with the same parameters as the channel in Figure 5.

The necessary condition we told about states that queue 1 can never behave like a processor sharing model, queue 2 can only if $x < 0.292893$, and queue 3 always can. This is, indeed, what we see on Figure 5. All plain lines are lower bounds and sometimes fit well simulation results. Queue 1 is "pushed" by 2 and then by 3; queue 2 is on its own (i.e., behaves like a processor sharing) until approximately $x = 0.2$, then is pushed by 3; queue 3 is always on its own. We see the necessary condition $x < 0.292893$ holds, but is not very tight. We hope this kind of study can furnish more detailed results and holds for other policies.

## 3. Localized empirical generator, entropy and the rate function.

3.1. *Localized empirical generator.* Take $x \in \mathbb{R}_+^{\mathcal{R}}$ and $D \in \mathbb{R}^{\Lambda(x)}$. We are interested in computing large deviations bounds of the form (1.3) (i.e., linear bounds as presented in Figure 3). In order to prove Theorem 1.2 we introduce a functional which allows one to measure how the different arrival rates should be modified in order that the rescaled process $Q_x^n$ follows a prescribed drift $D$. Moreover, the explanation exposed just after the statement of Theorem 1.2 suggests that the transition rates of routes indexed by $\Lambda_2(x)$ should not be modified and so it is useless to measure the arrivals on routes belonging to $\Lambda_2(x)$. Let us introduce the localized empirical generator at point $x$, as well as suitable state spaces associated to this process:

DEFINITION 3.1 (Localized empirical generators). Let $\Lambda$ be a face and denote:

- $A_{ij}(t)$, the number of arrivals on route $ij$ till $t$;
- the restriction $A^\Lambda(t) \stackrel{\text{def}}{=} (A_{ij}(t), \; ij \in \Lambda \cup \Lambda_1)$;
- $G_t^\Lambda = (\frac{1}{t} A^\Lambda(t), \frac{Q_t - Q_0}{t})$, the localized empirical generator on the face $\Lambda$.

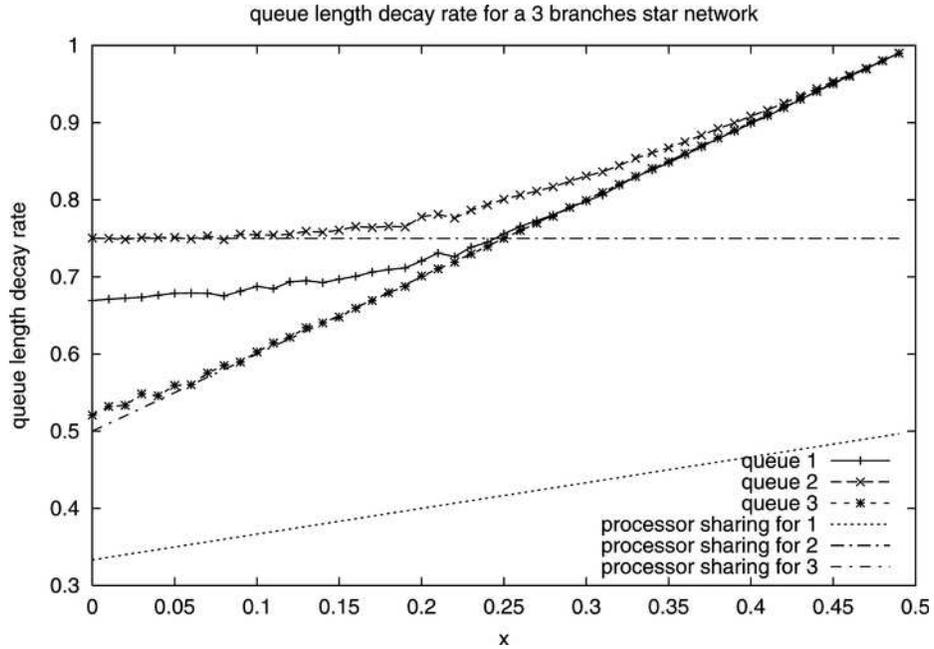

FIG. 5. *Comparison between decay rates obtained by simulation to processor sharing models with the equivalent parameters. Simulations were stopped at time $T = 10^8$.*

The set $\Gamma^\Lambda$ of localized empirical generators is the set of elements $(A^\Lambda, D)$ with $D \in \mathbb{R}^\mathcal{R}$ satisfying:

(3.1) $\quad$ (i) $\quad a_{ij} \geq 0 \; \forall ij \in \Lambda \cup \Lambda_1,$
$\qquad\quad$ (ii) $\quad a_{ij} - D_{ij} \geq 0 \; \forall ij \in \Lambda \cup \Lambda_1.$

The space $\Gamma^\Lambda$ is equipped with the distance $d$ defined by

$$d(G, G') \stackrel{\text{def}}{=} \sum_{ij \in \Lambda \cup \Lambda_1} |a_{ij} - a'_{ij}| + \sum_{ij \in \mathcal{R}} |D_{ij} - D'_{ij}| \qquad \forall G, G' \in \Gamma^\Lambda.$$

The inequalities (i) and (ii) in (3.1) refer, respectively, to the mean number of arrivals $a_{ij}$ and to the mean number of deconnections per unit time, $a_{ij} - D_{ij}$ being positive.

Since it is difficult to analyze at first the behavior of $Q(t)$ as in (1.3), we shall first establish large deviation bounds for the event

$$(3.2) \quad E^{(n)}_{\tau,\delta,y}(x, G) \stackrel{\text{def}}{=} \left\{ G^{\Lambda(x)}_{n\tau} \in B(G, \delta), \sup_{t \in [0,n\tau]} |Q(t, y) - nx - Dt| < \delta n \right\},$$

where $B(G, \delta)$ is the ball of center $G$ and radius $\delta$ [within the metric space $(\Gamma^{\Lambda(x)}, d)$]. As it will emerge, strong constraints must be imposed on $G$ in order that the event $E^{(n)}_{\tau,\delta,y}(x, G)$ occurs at a large deviation scale. More precisely, the arrivals must be cut on routes belonging to $\Lambda_1(x)$.

LEMMA 3.2.  *Take $x \in \mathbb{R}^\mathcal{R}_+$ and $G = (A, D) \in \Gamma^{\Lambda(x)}$, such that $D \in \mathbb{R}^{\Lambda(x)}$. If there exist $m$ and $p$ such that*

$$x_m = 0 \quad \text{and} \quad x_p > 0, \quad \text{and} \quad a_{pm} > 0,$$

*then $E^{(n)}_{\tau,\delta,y}(x, G)$ almost never occurs at a large deviation scale, that is,*

$$(3.3) \quad \lim_{\tau,\delta,\epsilon \to 0} \limsup_{n \to \infty} \frac{1}{n\tau} \sup_{|y - nx| < \epsilon n} \log \mathbb{P}[E^{(n)}_{\tau,\delta,y}(x, G)] = -\infty.$$

PROOF. The proof relies on a change of measure, as in Section 4.1. In fact, on $E^{(n)}_{\tau,\delta,y}(x, G)$ the service rate on route $pm$ tends to 0 when the different limits are taken. Since on $E^{(n)}_{\tau,\delta,y}(x, G)$, the arrival process is not cut on route $pm$, the cost to keep the component $pm$ of the rescaled process near 0 is infinite. Details are similar to the proof of the upper bound (see Section 4.2) and are omitted. $\square$

Lemma 3.2 states that in order to prove large deviation bounds for the localized empirical generator, it will be sufficient to deal with the following subspace of $\Gamma^{\Lambda(x)}$.

DEFINITION 3.3. $\mathcal{G}^\Lambda$ denotes the set of localized empirical generators $(A^\Lambda, D)$ such that:

(3.4)
- (i) $D \in \mathbb{R}^\Lambda$,
- (ii) $a_{ij} = 0 \; \forall ij \in \Lambda_1$,
- (iii) $a_{ij} > 0$ and $a_{ij} - D_{ij} > 0 \; \forall ij \in \Lambda$.

In this setting $\overline{\mathcal{G}}^\Lambda$ will represent the closure of $\mathcal{G}^\Lambda$.

Owing to Lemma 3.2, it is sufficient to deal with empirical generators satisfying (ii). In order to prove the large deviation local bounds, it will be sufficient to deal with empirical generators such that arrival and service rates are not cut, for $ij \in \Lambda(x)$, hence, condition (iii). A simple continuity argument will allow to extend the bounds obtained for $G \in \mathcal{G}^{\Lambda(x)}$ to $G \in \overline{\mathcal{G}}^{\Lambda(x)}$.

3.2. *Correspondence between localized empirical generators and star networks.* Let $G = (A, D) \in \overline{\mathcal{G}}^\Lambda$ be a localized empirical generator. It is associated with a unique localized star network $(\tilde{\lambda}_{ij}, \tilde{\mu}_{ij}(y), y \in \mathbb{R}_+^\mathcal{R})$ by the following relations:

(3.5)
- (i) $\tilde{\lambda}_{ij} = a_{ij} \; \forall ij \in \Lambda \cup \Lambda_1$,
- (ii) $\tilde{\lambda}_{ij} = \lambda_{ij} \; \forall ij \in \Lambda_2$,
- (iii) $\tilde{\mu}_{ij}(y) = \tilde{\lambda}_{ij} - D_{ij} \; \forall ij \in \Lambda \cup \Lambda_1, \; \forall y \in \mathbb{R}_+^\mathcal{R}$,
- (iv) $\tilde{\mu}_{ij}(y) = \mu_{ij}(y) \; \forall ij \in \Lambda_2, \; \forall y \in \mathbb{R}_+^\mathcal{R}$.

Let us describe the behavior of this network when it starts from $x$ [with $\Lambda = \Lambda(x)$]. In this case the routes belonging to $\Lambda_2$ behave as a star network of the type presently studied and the parameters of the routes belonging to this set are left unchanged. Moreover, they are independent from the rest of the network. Indeed, if $ij \in \Lambda_2$, then $x_{ik} = 0$ for all $k$ such that $ik \notin \Lambda_2$ (actually $ik \in \Lambda_1$; see Figure 2), hence, the constraints imposed on $G$ insures that $\tilde{\lambda}_{ik} = 0$. Hence,

$$\mu_{ij}(Q(s)) = Q_{ij}(s) \frac{\mu_i}{\sum_{ik \in \Lambda_2^c} Q_{ik}(s)} \wedge \frac{\mu_j}{\sum_{jk \in \Lambda_2^c} Q_{jk}(s)} \qquad \forall ij \in \Lambda_2^c(x),$$

proving the asserted independence. Moreover, the network consisting of the routes belonging to $\Lambda_2$ is ergodic when the initial network is. Indeed, for all ergodic channel $i$ (see Figure 2),

$$\sum_{j: ij \in \Lambda_2} \frac{\tilde{\lambda}_{ij}}{\tilde{\mu}_{ij}} = \sum_{j: ij \in \Lambda_2} \frac{\lambda_{ij}}{\mu_{ij}} \leq \sum_{j \in \mathcal{S}} \frac{\lambda_{ij}}{\mu_i} < C_i.$$

Besides, routes belonging to $\Lambda$ behave like independent $M/M/1$ queues up to the initial conditions, whereas the routes indexed by $\Lambda_1$ remain null. Now, the parameters have been chosen so that:

LEMMA 3.4. *Assume that $Q$ is ergodic. Let $x \in \mathbb{R}_+^{\mathcal{R}}$, $G = (A, D) \in \mathcal{G}^{\Lambda(x)}$, a localized empirical generator, and denote $\widetilde{\mathbb{P}}$ the law of its associated star network. Then, for all $\tau$,*

$$\lim_{\delta,\epsilon \to 0} \inf_{|y-nx|<\epsilon n} \liminf_{n \to \infty} \widetilde{\mathbb{P}}[E_{\tau,y}^{(n)}(x,G) \cap \{A_{ij}(n\tau) = 0 \ \forall ij \in \Lambda_1(x)\}] = 1.$$

PROOF. The proof is omitted: it is a classical fluid limit. □

### 3.3. *Entropy.*

DEFINITION 3.5. Let $x \in \mathbb{R}_+^{\mathcal{R}}$, $R(x) = (\lambda_{ij}, \mu_{ij}(x))$ denotes the generator of the star network at $x$, $G = (A, D) \in \mathcal{G}^{\Lambda(x)}$ be a localized generator and $(\tilde{\lambda}_{ij}, \tilde{\mu}_{ij}(y), y \in \mathbb{R}_+^{\mathcal{R}})$ its representation as a star network. The relative entropy of $G$ with respect to $R(x)$ is

$$H(G\|R(x)) = \sum_{ij \in \Lambda(x) \cup \Lambda_1(x)} I_p(\tilde{\lambda}_{ij}\|\lambda_{ij}) + I_p(\tilde{\mu}_{ij}\|\mu_{ij}(x)),$$

where $I_p(\nu\|\lambda)$ is the relative entropy of Poisson processes of intensities $\nu$ and $\lambda$ defined by

(3.6) $$I_p(\nu\|\lambda) \stackrel{\text{def}}{=} \nu \log \frac{\nu}{\lambda} - \nu + \lambda,$$

with the convention $\frac{0}{0} = 0$ and $0 \log 0 = 0$.

The entropy has an easy interpretation in terms of information theory: it can be defined as the *mean information gain*. $H(\cdot\|R)$ is decomposed as the sum of the information gain for the arrivals $I_p(\tilde{\lambda}_{ij}\|\lambda_{ij})$, the information gain for the service time $I_p(\tilde{\mu}_{ij}\|\mu_{ij}(x))$.

LEMMA 3.6. *For fixed $x$, $H(\cdot\|R(x))$ is continuous on $\overline{\mathcal{G}}^{\Lambda(x)}$.*

PROOF. It is an easy consequence of the expression (3.6). □

### 3.4. *The local rate function $L(x,D)$.*

DEFINITION 3.7. The local rate function $L(x,D)$ is defined by

(3.7) $$L(x,D) \stackrel{\text{def}}{=} \inf_{G \in f_{\Lambda(x)}^{-1}(D)} H(G\|R(x)) \qquad \forall D \in \mathbb{R}^{\Lambda(x)},$$

where $f_{\Lambda(x)} : \mathcal{G}^{\Lambda(x)} \mapsto \mathbb{R}^{\Lambda(x)}$ is the projection $f_{\Lambda(x)}(G) = D$.

It appears that $L(x, D)$ is the cost for a set of $M/M/1$ independent queues indexed by $\Lambda(x) \cup \Lambda_1(x)$ to follow the prescribed drift $D$ when the queues are far from all boundaries. A simple computation yields

$$l(D\|\lambda,\mu) \stackrel{\text{def}}{=} D\log\left(\frac{D+\sqrt{D^2+4\lambda\mu}}{2\lambda}\right) + \lambda + \mu - \sqrt{D^2+4\lambda\mu} \geq 0$$

for the cost that a $M/M/1$ queue with parameters $\lambda$ and $\mu$ follows the drift $D$ [see, e.g., Shwartz and Weiss ((1995))]. Using this remark and the identity $l(0\|\lambda,0) = \lambda$, one can deduce the explicit representation (1.4) for $L(x,D)$ [which is equal to (3.8) under the constraint $D_{ij} = \mu_{ij}(x) = 0$ for $ij \in \Lambda_1(x)$].

In equations (1.4) and (3.7), $L(x,D)$ is only defined for $D \in \mathbb{R}^{\Lambda(x)}$. In order to study the properties of the rate function $I_T(\cdot)$, it is convenient to extend the definition of $L(x,D)$ for all $D$, such that $D_{ij} \geq 0$ for all $ij \in \Lambda^c(x)$ by

$$(3.8) \qquad L(x,D) \stackrel{\text{def}}{=} \sum_{ij \in \Lambda(x) \cup \Lambda_1(x)} l(D_{ij}\|\lambda_{ij}, \mu_{ij}(x)).$$

PROPOSITION 3.8. *The local rate function $L(x,D)$ possesses the following properties.*

(i) *It is positive, finite, strictly convex and continuous with respect to $D$, such that $D_{ij} \geq 0$ for all $ij \in \Lambda(x)$. It has compact level sets;*

(ii) *there exists $M \in \mathbb{R}$ such that*

$$L(x,D) \geq \tfrac{1}{2}\|D\|\log\|D\| \qquad \forall\, x \in \mathbb{R}_+^{\mathcal{R}},\ \forall\, \|D\| \geq M;$$

(iii) *for a fixed $D$ and a prescribed face $\Lambda$, $L(x,D)$ is continuous for $x \in \Lambda$ [see (1.2)];*

(iv) *$L(x,D)$ is jointly lower semicontinuous w.r.t. $x$ and $D$.*

PROOF. Properties (i) and (ii) are obvious from (3.8).

(iii) is clear from (3.8), noting that the functions $\mu_{ij}(x)$, $ij \in \Lambda$, are continuous for $x$ belonging to the face $\Lambda$. Moreover, $\Lambda_1(x) = \Lambda_1$ is constant for $x \in \Lambda$.

Let $(x^{(n)}, D^{(n)})$ tends to $(x,D)$. First, it is clear that for $n$ large enough, $\Lambda(x) \subset \Lambda(x^{(n)})$ and also $\Lambda(x) \cup \Lambda_1(x) \subset \Lambda(x^{(n)}) \cup \Lambda_1(x^{(n)})$. Hence, since $l$ is positive, for sufficiently large $n$,

$$(3.9) \qquad L(x^{(n)}, D^{(n)}) \geq \sum_{ij \in \Lambda(x) \cup \Lambda_1(x)} l(D_{ij}^{(n)}\|\lambda_{ij}, \mu_{ij}(x^{(n)})).$$

Now, $\lambda_{ij} > 0$ (since $ij \in \mathcal{R}$) so that $l(\cdot\|\lambda_{ij},\cdot)$ is continuous. Moreover, $\mu_{ij}(x^{(n)}) \to \mu_{ij}(x)$ $\forall\, ij \in \Lambda(x) \cup \Lambda_1(x)$. Therefore the right-hand side of (3.9) converges to $L(x,D)$ and the lower semicontinuity (iv) is proved. $\square$

3.5. *The sample path rate function $I_T(\cdot)$.*

PROPOSITION 3.9. *The rate function $I_T(\cdot)$ defined in (1.6) possesses the following properties:*

(i) *Assume $I_T(\varphi) \leq K$ for some $K$. Then, for all $\epsilon > 0$, there exists $\delta > 0$ independent of $\varphi$, such that for any collection of nonoverlapping intervals $[t_j, t_{j+1}]$ in $[0,T]$ with $\sum_j t_{j+1} - t_j = \delta$,*

$$\sum_j |\varphi(t_{j+1}) - \varphi(t_j)| \leq \epsilon;$$

(ii) *$I_T(\cdot)$ is lower semicontinuous in $(D([0,T], \mathbb{R}_+^{\mathcal{R}}), d_d)$;*

(iii) *for $C \subset \mathbb{R}_+^{\mathcal{R}}$ compact, $\bigcup_{x \in C} \Phi_x(K)$ is compact in $\mathcal{C}([0,T], \mathbb{R}_+^{\mathcal{R}})$ [see (1.7) for the definition of the level set $\Phi_x(K)$];*

(iv) *consider an absolutely continuous function $\varphi$ with $I_T(\varphi) < \infty$. Then, for all $\epsilon > 0$, there exists a piecewise linear function $\varphi_\epsilon$ such that:*

(a) $d_c(\varphi_\epsilon, \varphi) \leq \epsilon$,
(b) $I_T(\varphi_\epsilon) \leq I_T(\varphi) + \epsilon$.

PROOF. One proves (i) using Proposition 3.8(ii) in a way similar to Lemma 5.18 of Shwartz and Weiss (1995).

In order to prove the lower semicontinuity of $I_T(\cdot)$, (i) shows it is sufficient to consider sequences of absolutely continuous functions. Since on $\mathcal{C}([0,T], \mathbb{R}_+^{\mathcal{R}})$, the metrics $d_c$ and $d_d$ are equivalent, one can use $d_c$. Now, using Proposition 3.8(ii), the fact that $L(x,D)$ is lower semicontinuous in $(x,D)$ and convex with respect to $D$ by Proposition 3.8, (ii) is proved by means of Theorem 3 of Section 9.1.4 in Ioffe and Tihomirov (1979).

(iii) is a consequence of (i) and (ii) [see Proposition 5.46 of Shwartz and Weiss (1995)].

The proof of (iv) is a simple adaptation of Proposition 6.3(4) of Delcoigne and de La Fortelle (2002). □

**4. Large deviations bounds for the localized empirical generator.** In this section, we aim at proving the following theorem:

THEOREM 4.1. *Let $x \in \mathbb{R}_+^{\mathcal{R}}$ and $G = (A, D) \in \overline{\mathcal{G}}^{\Lambda(x)}$ be a localized generator. Then*

$$-H(G\|R(x)) = \lim_{\tau,\delta,\epsilon \to 0} \inf_{|y-nx|<\epsilon n} \liminf_{n \to \infty} \frac{1}{n\tau} \log \mathbb{P}[E_{\tau,\delta,y}^{(n)}(x,G)]$$

$$= \lim_{\tau,\delta,\epsilon \to 0} \sup_{|y-nx|<\epsilon n} \limsup_{n \to \infty} \frac{1}{n\tau} \log \mathbb{P}[E_{\tau,\delta,y}^{(n)}(x,G)],$$

where $E^{(n)}_{\tau,\delta,y}(x,G)$ is the event defined in (3.2). Moreover, if a face $\Lambda$ and a drift $D \in \mathbb{R}^\Lambda$ are fixed, then the preceding limit in $\tau$ is uniform with respect to $x$ in compact sets of $\Lambda$ (see Definition 1.1).

4.1. *An exponential change of measure.* Fix an empirical generator $G = (A, D) \equiv (\tilde{\lambda}_{ij}, \tilde{\mu}_{ij}(y), y \in \mathbb{R}^{\mathcal{R}}_+) \in \overline{\mathcal{G}}^{\Lambda(x)}$ and denote by the following:

- $N_t$, the number of jumps of the process till $t$.
- $Q(k) = \{Q_{ij}(k), \ i,j \in \mathcal{S}\}$, the embedded Markov chain at time $k \in \mathbb{N}$. We shall distinguish between discrete and continuous time by using $k$ for discrete and $s$ or $t$ for continuous time.

Define the following:

- The mapping $h: \mathbb{Z}^{\mathcal{R}}_+ \times \mathbb{Z}^{\mathcal{R}}_+ \mapsto \mathbb{R}$ by

$$h(x,y) \stackrel{\text{def}}{=} \begin{cases} \log \dfrac{\tilde{\lambda}_{ij}}{\lambda_{ij}}, & \text{if } y - x = e_{ij} \text{ and } \tilde{\lambda}_{ij} > 0, \\ \log \dfrac{\tilde{\mu}_{ij}(x)}{\mu_{ij}(x)}, & \text{if } y - x = -e_{ij} \text{ and } \tilde{\mu}_{ij}(x) > 0, \\ 0, & \text{otherwise}. \end{cases}$$

- The compensator $K: \mathbb{Z}^{\mathcal{R}}_+ \mapsto \mathbb{R}$ by

(4.1)
$$K(x) \stackrel{\text{def}}{=} \sum_{y \in \mathbb{Z}^{\mathcal{R}}_+} q(x,y)(e^{h(x,y)} - 1)$$
$$= \sum_{ij \in \mathcal{R}} (\tilde{\lambda}_{ij} - \lambda_{ij}) + \sum_{ij \in \mathcal{R}} (\tilde{\mu}_{ij}(x) - \mu_{ij}(x)).$$

- The process

$$\mathcal{M}_t \stackrel{\text{def}}{=} \exp\left\{ \sum_{k=0}^{N_t - 1} h(Q(k), Q(k+1)) - \int_0^t K(Q(s))\, dv \right\}.$$

Note that the compensator is always bounded, so that $\mathcal{M}_t$ takes only finite values. Since $K$ has been exactly defined so that

$$K(x) = \frac{d}{dt}\mathbb{E}\left[\exp\left\{ \sum_{k=0}^{N_t-1} h(Q(k,x), Q(k+1,x)) \right\}\right]_{t=0},$$

it is easily checked that the derivative of $\mathbb{E}[\mathcal{M}_t]$ at $t = 0$ is null (note that the derivative is independent of $\Lambda$, so that it is dropped). Then using the Markov property, one can get that the derivative is null for all $t \geq 0$, so that

$$\mathbb{E}[\mathcal{M}_t] = 1.$$

Using again the Markov property, this proves that

$$\mathbb{E}[\mathcal{M}_t | \mathcal{F}_s] = \mathcal{M}_s \qquad \text{for all } t \geq s \geq 0,$$

hence, $\{\mathcal{M}_t,\ t \geq 0\}$ is a martingale w.r.t. the natural filtration $\mathcal{F}_t$.

Then define a new probability measure by

$$\widetilde{\mathbb{P}}[B] \stackrel{\text{def}}{=} \mathbb{E}[\mathbb{1}_{\{B\}} \mathcal{M}_t] \qquad \forall\, B \in \mathcal{F}_t.$$

It is a matter of routine to show that under $\widetilde{\mathbb{P}}$, $X$ is again a Markov process. In fact, under $\widetilde{\mathbb{P}}$, the system behaves like a star network, where the arrival and the service rates at node $ij$ are respectively given by $\tilde{\lambda}_{ij}$ and $\tilde{\mu}_{ij}(y)$ (whence the notation).

REMARK. The probability measure $\mathbb{P}$ is not necessarily absolutely continuous with respect to $\widetilde{\mathbb{P}}$. This is the case, for instance, if for some $ij \in \mathcal{R}$, $\tilde{\lambda}_{ij} = 0$ (whereas $\lambda_{ij} > 0$).

4.2. *Proof of the upper bound of Theorem* 4.1. Since $\mathbb{P}$ is not necessarily absolutely continuous with respect to $\widetilde{\mathbb{P}}$, in order to prove the upper bound, we introduce a sequence of change of measure $\{\widetilde{\mathbb{P}}^{(\eta)},\ \eta > 0\}$ such that

$$\begin{aligned}
\tilde{\lambda}_{ij}^{(\eta)} > 0 \quad &\text{and} \quad \lim_{\eta \to 0} \tilde{\lambda}_{ij}^{(\eta)} = \tilde{\lambda}_{ij} \qquad \forall\, ij \in \Lambda(x) \cup \Lambda_1(x), \\
\tilde{\mu}_{ij}^{(\eta)} > 0 \quad &\text{and} \quad \lim_{\eta \to 0} \tilde{\mu}_{ij}^{(\eta)} = \tilde{\mu}_{ij}(x) \qquad \forall\, ij \in \Lambda(x).
\end{aligned}$$

In this setting, $\{\mathcal{M}_t^{(\eta)},\ t \geq 0\}$ is the martingale defining $\widetilde{\mathbb{P}}^{(\eta)}$ with respect to $\mathbb{P}$, and $h^{(\eta)}(x,y)$ and $K^{(\eta)}(x)$ are the functions used to defined $\mathcal{M}_t^{(\eta)}$ according to Section 4.1. Now, $\widetilde{\mathbb{P}}^{(\eta)}$ and $\mathbb{P}$ are mutually absolutely continuous and

$$(4.2) \qquad \mathbb{P}[E_{\tau,\delta,y}^{(n)}(x,G)] = \widetilde{\mathbb{E}}^{(\eta)}\left[\mathbb{1}_{\{E_{\tau,\delta,y}^{(n)}(x,G)\}} (\mathcal{M}_{n\tau}^{(\eta)})^{-1}\right].$$

Let us majorize $(\mathcal{M}_{n\tau}^{(\eta)})^{-1}$ on $E_{\tau,\delta,y}^{(n)}(x,G)$ when $|y - nx| < \epsilon n$. First, recalling

$$\tilde{\lambda}_{ij} = \lambda_{ij} \qquad \text{for } ij \in \Lambda_2$$

and

$$\tilde{\mu}_{ij}(y) = \mu_{ij}(y) \qquad \text{for } ij \in \Lambda_1 \cup \Lambda_2 \text{ and } y \in \mathbb{R}_+^{\mathcal{R}}$$

one has the following bounds:

$$-\sum_{k=0}^{N_{n\tau}-1} h^{(\eta)}(Q(k), Q(k+1))$$
$$\leq -n\tau \left( \sum_{ij \in \Lambda(x)} \tilde{\mu}_{ij} \log \frac{\tilde{\mu}_{ij}^{(\eta)}}{\sup_{s \in [0,n\tau]} \mu_{ij}(Q(s))} \right.$$
$$\left. + \sum_{ij \in \Lambda(x) \cup \Lambda_1(x)} \tilde{\lambda}_{ij} \log \frac{\tilde{\lambda}_{ij}^{(\eta)}}{\lambda_{ij}} \right)$$
$$+ n\tau \delta \left( \sum_{ij \in \Lambda(x)} \left| \log \frac{\tilde{\mu}_{ij}^{(\eta)}}{\inf_{s \in [0,n\tau]} \mu_{ij}(Q(s))} \right| + \sum_{ij \in \Lambda(x) \cup \Lambda_1(x)} \left| \log \frac{\tilde{\lambda}_{ij}^{(\eta)}}{\lambda_{ij}} \right| \right).$$

(4.3)

Moreover, the compensator $K$ is bounded in (4.1) by

(4.4)
$$\int_0^{n\tau} K^{(\eta)}(Q(s))\, ds$$
$$\leq n\tau \sum_{ij \in \Lambda(x) \cup \Lambda_1(x)} (\tilde{\lambda}_{ij}^{(\eta)} - \lambda_{ij}) + n\tau \sum_{ij \in \Lambda(x)} \left( \tilde{\mu}_{ij}^{(\eta)} - \inf_{s \in [0,n\tau]} \mu_{ij}(Q(s)) \right).$$

Besides, on $E_{\tau,\delta,y}^{(n)}(x, G)$, we have for $ij \in \Lambda(x)$

$$0 < \mu_{ij}(x) = \lim_{\tau,\delta,\epsilon \to 0} \inf_{|y-nx|<\epsilon n} \liminf_{n \to \infty} \inf_{s \in [0,n\tau]} \mu_{ij}(Q(s,y))$$
$$= \lim_{\tau,\delta,\epsilon \to 0} \sup_{|y-nx|<\epsilon n} \limsup_{n \to \infty} \sup_{s \in [0,n\tau]} \mu_{ij}(Q(s,y)).$$

Finally, majorizing $\mathbb{1}_{\{E_{\tau,\delta,y}^{(n)}(x,G)\}}$ by 1, bounding $\mathcal{M}_{n\tau}^{(\eta)}$ using (4.3), (4.4) and (4.5) and taking into account the order in which the different limits are taken, the representation formula (4.2) yields

$$\lim_{\tau,\delta,\epsilon \to 0} \limsup_{n \to \infty} \frac{1}{n\tau} \sup_{|y-nx|<\epsilon n} \log \mathbb{P}[E_{\tau,\delta,y}^{(n)}(x, G)]$$
$$\leq - \sum_{ij \in \Lambda(x) \cup \Lambda_1(x)} \tilde{\lambda}_{ij} \log \frac{\tilde{\lambda}_{ij}^{(\eta)}}{\lambda_{ij}} - \tilde{\lambda}_{ij}^{(\eta)} + \lambda_{ij}$$
$$- \sum_{ij \in \Lambda(x)} \tilde{\mu}_{ij}^{(\eta)} \log \frac{\tilde{\mu}_{ij}^{(\eta)}}{\mu_{ij}(x)} - \tilde{\mu}_{ij}^{(\eta)} + \mu_{ij}(x).$$

The proof of the upper bound is concluded letting $\eta$ tend to 0.

4.3. *Proof of the lower bound of Theorem* 4.1. Take $G \in \mathcal{G}^{\Lambda(x)}$ and denote the event (appearing in Lemma 3.4)

$$F_{\tau,\delta,y}^{(n)}(x,G) \stackrel{\text{def}}{=} E_{\tau,y}^{(n)}(x,G) \cap \{A_{ij}(n\tau) = 0 \ \forall ij \in \Lambda_1(x)\}.$$

Although $\mathbb{P}$ is not absolutely continuous w.r.t. $\widetilde{\mathbb{P}}$, by definition of $\mathcal{G}^{\Lambda(x)}$, $\tilde{\lambda}_{ij} > 0$ and $\tilde{\mu}_{ij} > 0 \ \forall ij \in \Lambda(x)$ so that $\mathbb{P}$ is absolutely continuous w.r.t. $\widetilde{\mathbb{P}}$ on $F_{\tau,\delta,y}^{(n)}(x,D)$ and

$$\mathbb{P}[E_{\tau,\delta,y}^{(n)}(x,G)] \geq \mathbb{P}[F_{\tau,\delta,y}^{(n)}(x,G)]$$

$$\geq \inf_{\omega \in F_{\tau,\delta,y}^{(n)}(x,D)} \mathcal{M}_{n\tau}^{-1}(\omega) \widetilde{\mathbb{P}}[F_{\tau,\delta,y}^{(n)}(x,G)].$$

By Lemma 3.4, $\widetilde{\mathbb{P}}[F_{\tau,\delta,y}^{(n)}(x,G)]$ tends to 1. Therefore, reversing the inequalities obtained for the upper bound yields

$$\lim_{\tau,\delta,\epsilon \to 0} \liminf_{n \to \infty} \frac{1}{n\tau} \inf_{|y-nx|<\epsilon n} \log \mathbb{P}[E_{\tau,\delta,y}^{(n)}(x,G)]$$

$$\geq - \sum_{ij \in \Lambda(x) \cup \Lambda_1(x)} \tilde{\lambda}_{ij} \log \frac{\tilde{\lambda}_{ij}}{\lambda_{ij}} - \tilde{\lambda}_{ij} + \lambda_{ij}$$

$$- \sum_{ij \in \Lambda(x)} \tilde{\mu}_{ij} \log \frac{\tilde{\mu}_{ij}}{\mu_{ij}(x)} - \tilde{\mu}_{ij} + \mu_{ij}(x).$$

This concludes the proof of the lower bound when $G \in \mathcal{G}^{\Lambda(x)}$.

Consider $G \in \overline{\mathcal{G}}^{\Lambda(x)}$ and define $G^{(\varepsilon)}$ by $\tilde{\lambda}_{ij}^{(\varepsilon)} \stackrel{\text{def}}{=} \tilde{\lambda}_{ij} + \varepsilon$ and $\tilde{\mu}_{ij}^{(\varepsilon)} \stackrel{\text{def}}{=} \tilde{\mu}_{ij} + \varepsilon$, for $ij \in \Lambda(x)$; otherwise, the coefficients are the same. Then $G^{(\varepsilon)}$ belongs to $\mathcal{G}^{\Lambda(x)}$ for $\varepsilon > 0$, it converges to $G$ and its entropy converges to $H(G\|R(x))$ by Lemma 3.6. Moreover, the drifts $D = (\tilde{\lambda}_{ij} - \tilde{\mu}_{ij})$ and $D^{(\varepsilon)} = (\tilde{\lambda}_{ij}^{(\varepsilon)} - \tilde{\mu}_{ij}^{(\varepsilon)})$ are equal.

For any $\varepsilon_0 > 0$ there exists $\varepsilon_1 > 0$ and $\delta_1 > 0$ such that, for all $0 < \varepsilon' < \varepsilon_1$ and $0 < \delta' < \delta_1$, $B(G^{(\varepsilon)}, \delta') \subset B(G, \delta)$ and $H(G^{(\varepsilon)}\|R(x)) \leq H(G\|R(x)) + \varepsilon_0$. For the sake of simplicity, we shall denote $G^{(\varepsilon)}$ by $G'$. Since $D' = D$, we get the *time uniform* inclusion

$$E_{\tau,\delta',y}^{(n)}(x,G') \subset E_{\tau,\delta,y}^{(n)}(x,G) \qquad \forall \tau \geq 0.$$

It yields, using the decrease of $E_{\tau,\delta',y}^{(n)}(x,G')$ with $\delta'$,

$$\lim_{\epsilon \to 0} \liminf_{n \to \infty} \frac{1}{n\tau} \inf_{|y-nx|<\epsilon n} \log \mathbb{P}[E_{\tau,\delta,y}^{(n)}(x,G)]$$

$$\geq \lim_{\delta',\epsilon \to 0} \liminf_{n \to \infty} \frac{1}{n\tau} \inf_{|y-nx|<\epsilon n} \log \mathbb{P}[E_{\tau,\delta',y}^{(n)}(x,G')] \qquad \forall \tau > 0.$$

Using the lower bound for $G' \in \mathcal{G}^{\Lambda(x)}$ and the uniformity over time of the previous bound, by letting $\tau$ tend to 0 we deduce that (depending on $\varepsilon_0$) there exists $\tau_0$ such that, for all $0 < \tau < \tau_0$,

$$\lim_{\epsilon \to 0} \liminf_{n \to \infty} \frac{1}{n\tau} \inf_{|y-nx|<\epsilon n} \log \mathbb{P}[E^{(n)}_{\tau,\delta,y}(x,G)] \geq -H(G'\|R(x)) - \varepsilon_0.$$

Now recall that the entropy is bounded (by continuity) when $G' \to G$ so that there is no problem when $\delta$ decreases, for all $0 < \tau < \tau_0$,

$$\lim_{\delta,\epsilon \to 0} \liminf_{n \to \infty} \frac{1}{n\tau} \inf_{|y-nx|<\epsilon n} \log \mathbb{P}[E^{(n)}_{\tau,\delta,y}(x,G)] \geq -H(G\|R(x)) - 2\varepsilon_0.$$

Since this is true for any $\varepsilon_0$, we get the lower bound for $G$,

$$\lim_{\tau,\delta,\epsilon \to 0} \liminf_{n \to \infty} \frac{1}{n\tau} \inf_{|y-nx|<\epsilon n} \log \mathbb{P}[E^{(n)}_{\tau,\delta,y}(x,G)] \geq -H(G'\|R(x)).$$

Theorem 4.1 is proved for any $G \in \overline{\mathcal{G}}^{\Lambda(x)}$.

The uniformity of the limit stated in Theorem 4.1 is easily checked. Nonetheless, this uniformity is clear as far as $x$ evolves on compact sets of some face $\Lambda$. Indeed, if $x_{ij}$ goes to 0 for some $ij \in \Lambda$, then $\mu_{ij}(x)$ possibly vanishes and difficulties can appear.

PROOF OF THEOREM 1.2. Now Theorem 4.1 implies the large deviations local bounds of Theorem 1.2. Moreover, if a face $\Lambda$ and a drift $D \in \mathbb{R}^\Lambda$ are fixed, then the limits in (1.3) in $\tau$ are uniform w.r.t. $x$ in compact sets of $\Lambda$. The proof relies on a simple adaptation of the contraction principle, similarly to the proof of Theorem 7.2 of Delcoigne and de La Fortelle (2002). Details are omitted. $\square$

**5. Sample path LDP.** The proof of the sample path LDP is done in two steps which are briefly recalled. Using Markov property, Theorem 1.2 and the continuity of $L(x,D)$ with respect to $x \in \Lambda(D)$ for fixed $D$, large deviations bounds are established for the probability that the process stays near some linear path.

PROPOSITION 5.1 (Linear bounds). *Let $x \in \mathbb{R}^\mathcal{R}_+$ and $D \in \mathbb{R}^\mathcal{R}$, satisfying $x + DT \in \mathbb{R}^\mathcal{R}_+$. Denote $\varphi$ the function such that $\varphi(t) = x + Dt$ for all $t \in [0,T]$. Then*

$$-I_T(\varphi) = \lim_{\delta,\epsilon \to 0} \liminf_{n \to \infty} \frac{1}{n} \inf_{|y-nx|<\epsilon n} \log \mathbb{P}\left[\sup_{t \in [0,T]} |Q(t,y) - n\varphi(t)| < \delta n\right]$$

$$= \lim_{\delta,\epsilon \to 0} \limsup_{n \to \infty} \frac{1}{n} \sup_{|y-nx|<\epsilon n} \log \mathbb{P}\left[\sup_{t \in [0,T]} |Q(t,y) - n\varphi(t)| < \delta n\right].$$

PROOF. Due to the fact that the intensity $\mu_{ij}(x)$ is not bounded away from 0, this proof is quite involved. This is where the technical uniform reachability condition of Dupuis and Ellis ((1995)) is used; but it does not hold in the present model, even if the final result is the same. It is discussed in some detail in Section 4 and Appendix B of Delcoigne and de La Fortelle (2001). □

*From linear paths to LDP.* The sample path local bounds of Theorem 1.3 are now proved for linear paths (Proposition 5.1). There are some steps to reach the LDP, which we outline here.

First, the local bounds are extended to piecewise linear paths. Using the Markov property, the proof looks very much like that of Proposition 5.1.

Second, the local bounds are extended to absolutely continuous paths with finite entropy, using the properties of $I_T(\cdot)$. Notably points (ii) and (iv) of Proposition 3.9 imply that for an absolutely continuous $\varphi$ with $I_T(\varphi) < \infty$, there exists a sequence $\{\varphi_n,\ n \geq 1\}$ of piecewise linear paths satisfying

$$\lim_{n\to\infty} d_c(\varphi_n, \varphi) = 0 \quad \text{and} \quad \lim_{n\to\infty} I_T(\varphi_n) = I_T(\varphi).$$

The next step is to prove the exponential tightness of the sequence $\{n^{-1}Q(nt, [nx]),\ n \geq 1\}$ over finite interval of time (uniformly for $x$ belonging to a compact set). This is done, for instance, in Dupuis, Ellis and Weiss (1991). Finally, Theorem 1.3 is proved. These last two steps use various properties of the rate function $I_T(\cdot)$ and Proposition 3.9. The reader is referred to Section 5 of Dupuis and Ellis (1995) for details.

**6. LDP without ergodicity assumption.** Theorem 1.2 states large deviation bounds for ergodic networks. However, at the expense of cumbersome notation, it is possible to compute these bounds directly without ergodicity assumption introducing a more detailed empirical generator. For the ease of the exposition, the study was first performed for ergodic systems. We show now how one can compute, in general, $L(x, D)$. The discussion after Theorem 1.2 explains why the main difficulty to overcome is to compute the cost for an arbitrary star network under the min policy to stay in a neighborhood of 0.

PROPOSITION 6.1. *Let $Q$ be not necessarily ergodic. For all $\tau \geq 0$,*

$$\lim_{\delta,\epsilon\to 0} \inf_{|y|<\epsilon n} \liminf_{n\to\infty} \frac{1}{n\tau} \log \mathbb{P}\left[\sup_{t\in[0,n\tau]} |Q(t,y)| < \delta n\right]$$

$$= \lim_{\delta,\epsilon\to 0} \sup_{|y|<\epsilon n} \limsup_{n\to\infty} \frac{1}{n\tau} \log \mathbb{P}\left[\sup_{t\in[0,n\tau]} |Q(t,y)| < \delta n\right].$$

*The common value of these limits is denoted by* $-L(0,0)$ *and*

$$(6.1) \quad L(0,0) = \inf_{\nu \in V} \sum_{ij \in \mathcal{R}} (\sqrt{P\lambda_{ij}} - \sqrt{\mu_{ij}\nu_{ij}})^2 = \inf_{\nu \in V} \sum_{ij \in \mathcal{R}} l(0\|\lambda_{ij}, \mu_{ij}\nu_{ij}),$$

*where* $l(\cdot\|\cdot,\cdot)$ *is defined in* (1.5) *and the set* $V$ *by*

$$(6.2) \qquad V \stackrel{def}{=} \left\{ \nu \in \mathbb{R}_+^{\mathcal{R}} : \sum_{j \in \mathcal{S}} \nu_{ij} \leq C_i \; \forall i \in \mathcal{S} \right\}.$$

Note that Proposition 6.1 is a bit stronger than equality (1.3) of Theorem 1.2 applied to $x = D = 0$, since the time $\tau$ is not necessarily short. Besides, the rate function $L(0,0)$ is not explicit, but is an algorithmically fairly simple problem since it is a convex program w.r.t. $\sqrt{\nu_{ij}}$.

6.1. *Proof of Proposition* 6.1. As in the ergodic case, the proof relies on four steps: the introduction of a suitable empirical generator, the association of a star network to each empirical generator, the proof of large deviation bounds for empirical generator and finally, the proof of Proposition 6.1 using an adaptation of the contraction principle.

6.1.1. *Empirical generator.* This process is a bit different than the one defined in the ergodic case (see Definition 3.1). It takes into account the sole case $x = D = 0$, but in the transient case.

DEFINITION 6.2. The empirical generator $G_t$ is the functional defined by

$$G_t \stackrel{def}{=} \left( \frac{1}{t} A(t), \frac{1}{t} \int_0^t \nu(Q(s)) \, ds \right),$$

where $\nu(x) \stackrel{def}{=} (\nu_{ij}(x), \; i,j \in \mathcal{S})$. The set $\Gamma$ of empirical generators is $\mathbb{R}_+^{\mathcal{R}} \times V$; its elements will be denoted by $G = (A, \nu)$. It is equipped with the distance $d$ defined by

$$d(G, G') \stackrel{def}{=} \sum_{ij \in \mathcal{R}} |a_{ij} - a'_{ij}| + \sum_{ij \in \mathcal{R}} |\nu_{ij} - \nu'_{ij}| \qquad \forall G, G' \in \Gamma.$$

Large deviation bounds are established for the event [similarly to (3.2)]

$$(6.3) \qquad E_{\tau,\delta,y}^{(n)}(G) \stackrel{def}{=} \left\{ G_{n\tau} \in B(G,\delta), \sup_{t \in [0,n\tau]} |Q(t,y)| < \delta n \right\},$$

where $B(G, \delta)$ is the ball of center $G$ and radius $\delta$. Roughly speaking, when $\nu_{pm} = 0$ the service rate are cut on route $pm$ and so some constraints must be imposed on $A$. More precisely:

LEMMA 6.3. *Take $G = (A, \nu) \in \Gamma$. If there exist $m$ and $p$ such that*

$$\nu_{pm} = 0 \quad and \quad a_{pm} > 0,$$

*then $E^{(n)}_{\tau,\delta,y}(G)$ almost never occurs at a large deviation scale, that is,*

$$\lim_{\tau,\delta,\epsilon \to 0} \limsup_{n \to \infty} \frac{1}{n\tau} \sup_{|y| < \epsilon n} \log \mathbb{P}[E^{(n)}_{\tau,\delta,y}(G)] = -\infty.$$

PROOF. The proof is similar to the proof of Lemma 3.2. □

By Lemma 6.3, it is enough to deal with the following subspace of $\Gamma$.

DEFINITION 6.4. $\mathcal{G}$ denotes the set of empirical generator $(A, \nu)$ such that:

(i) $a_{ij} = 0$, when $\nu_{ij} = 0$,
(ii) $\sum_j \nu_{ij} < C_i \; \forall i$.

$\overline{\mathcal{G}}$ stands for the closure of $\mathcal{G}$.

6.1.2. *Correspondance between empirical generators and star networks.*
Let $G = (A, \nu) \in \overline{\mathcal{G}}$. It is associated arrival and departure rates

$$\tilde{\lambda}_{ij} \stackrel{\text{def}}{=} a_{ij} \qquad \forall ij \in \mathcal{R},$$

$$\tilde{\mu}_{ij}(y) \stackrel{\text{def}}{=} \tilde{\mu}_{ij} y_{ij} \frac{C_i}{y_i} \wedge \frac{C_j}{y_j} \mathbb{1}_{\{y_{ij} > 0\}} \qquad \forall ij \in \mathcal{R}, \; \forall y \in \mathbb{R}^{\mathcal{R}}_+,$$

where

$$\tilde{\mu}_{ij} \stackrel{\text{def}}{=} \begin{cases} \dfrac{\tilde{\lambda}_{ij}}{\nu_{ij}}, & \forall ij \text{ such that } \nu_{ij} > 0, \\ 0, & \text{otherwise.} \end{cases}$$

Then $(\tilde{\lambda}_{ij}, \tilde{\mu}_{ij}(y), \; y \in \mathbb{R}^{\mathcal{R}}_+)$ simply describes a star network under the min policy where the arrivals intensity and the duration of calls on route $ij$ are respectively given by $\tilde{\lambda}_{ij}$ and $\tilde{\mu}_{ij}$.

Similarly to Lemma 3.4, we now prove the following lemma:

LEMMA 6.5. *Let $G = (A, \nu) \in \mathcal{G}$ and $\widetilde{\mathbb{P}}$ the law of its associated star network. Then $Q$ is ergodic under $\widetilde{\mathbb{P}}$. Besides, for all $\tau$,*

(6.4) $$\lim_{\delta,\epsilon \to 0} \inf_{|y| < \epsilon n} \liminf_{n \to \infty} \widetilde{\mathbb{P}}[E^{(n)}_{\tau,\delta,y}] = 1.$$

PROOF. Since $G \in \mathcal{G}$, the ergodicity condition (1.1) are easily checked for $(\tilde{\lambda}_{ij}, \tilde{\mu}_{ij})$, so that $Q$ is ergodic under $\widetilde{\mathbb{P}}$. Moreover, a straight application of the ergodic theorem yields

$$(6.5) \qquad \lim_{t \to \infty} \frac{1}{t} \int_0^t \nu_{ij}(Q(s)) \, ds = \frac{\tilde{\lambda}_{ij}}{\tilde{\mu}_{ij}} = \nu_{ij} \qquad \forall ij.$$

Equation (6.4) is, thus, just a statement about fluid limits. □

### 6.1.3. Entropy and local bounds.

DEFINITION 6.6 (Entropy). Let $G = (A, D) \in \overline{\mathcal{G}}$ be an empirical generator and $(\tilde{\lambda}_{ij}, \tilde{\mu}_{ij})$ its representation as a star network. The relative entropy of $G$ with respect to $R$, the generator of the initial star network is

$$H(G\|R) = \sum_{ij \in \mathcal{R}} (I_p(\tilde{\lambda}_{ij}\|\lambda_{ij}) + I_p(\tilde{\lambda}_{ij}\|\nu_{ij}\mu_{ij})),$$

where $I_p$ is the entropy of Poisson processes defined in (3.6).

PROPOSITION 6.7. Let $G = (A, \nu) \in \overline{\mathcal{G}}$ be an empirical generator. Then

$$-H(G\|R) = \lim_{\delta,\epsilon \to 0} \inf_{|y| < \epsilon n} \liminf_{n \to \infty} \frac{1}{n\tau} \log \mathbb{P}[E^{(n)}_{\tau,\delta,y}(G)]$$

$$= \lim_{\delta,\epsilon \to 0} \sup_{|y| < \epsilon n} \limsup_{n \to \infty} \frac{1}{n\tau} \log \mathbb{P}[E^{(n)}_{\tau,\delta,y}(G)],$$

where $E^{(n)}_{\tau,\delta,y}(G)$ is the event defined in (6.3).

PROOF. The proof is similar to that of Theorem 4.1 and will not be repeated. Note simply that the lower bound is first proved for $G \in \mathcal{G}$ using, in particular, Lemma 6.5. It is then extended to all $G \in \overline{\mathcal{G}}$ using the continuity of the entropy $H$. □

PROOF OF PROPOSITION 6.1. Details are similar to the proof of Theorem 1.2 and, thus, omitted. Note that

$$L(0,0) = \inf_{G \in \overline{\mathcal{G}}} H(G\|R).$$

Taking $G = (A, \nu) \in \overline{\mathcal{G}}$ and minimizing w.r.t. $A$ yields (6.1). □

Taking into account Proposition 6.1, this leads to the following expression for $L(x, D)$ for a network without ergodicity condition and for $D \in \mathbb{R}^{\Lambda(x)}$:

$$(6.6) \quad L(x,D) = \sum_{ij \in \Lambda(x) \cup \Lambda_1(x)} l(D_{ij}\|\lambda_{ij}, \mu_{ij}(x)) + \inf_{\nu \in V} \sum_{ij \in \Lambda_2(x)} l(0\|\lambda_{ij}, \mu_{ij}\nu_{ij}),$$

where $V$ is defined in (6.2).

REMARK. At the expense of heavier notation, this theorem could have been derived at once, as in Section 4 studying the following more detailed empirical generator

$$L_t = \left(\frac{1}{t}A(t), \frac{1}{t}\int_0^t \nu_{\Lambda_2(x)}(Q(s))\,ds, \frac{1}{t}Q_t\right),$$

where $\nu_{\Lambda_2(x)} = (\nu_{ij}(y), ij \in \Lambda_2(x), y \in \mathbb{R}^{\mathcal{R}}_+)$.

## REFERENCES

ATAR, R. and DUPUIS, P. (1999). Large deviations and queueing networks: Methods for rate function identification. *Stochastic Process. Appl.* **84** 255–296.

BASKETT, F., CHANDY, K. M., MUNTZ, R. R. and PALACIOS, F. G. (1975). Open, closed, and mixed networks of queues with different classes of customers. *J. Assoc. Comput. Mach.* **22** 248–260. MR1719274

DELCOIGNE, F. and DE LA FORTELLE, A. (2001). Large deviations problems for star networks: The min policy. Technical Report 4143, INRIA.

DELCOIGNE, F. and DE LA FORTELLE, A. (2002). Large deviations rate function for polling systems. *Queueing Syst. Theory Appl.* **41** 13–44. MR1911125

DEMBO, A. and ZEITOUNI, O. (1998). *Large Deviations Techniques and Applications*, 2nd ed. Springer, New York. MR1619036

DEUSCHEL, J.-D. and STROOCK, D. W. (1989). *Large Deviations*. Academic Press, Boston. MR997938

DUPUIS, P. and ELLIS, R. S. (1992). Large deviations for Markov processes with discontinuous statistics. II. Random walks. *Probab. Theory Related Fields* **91** 153–194. MR1147614

DUPUIS, P. and ELLIS, R. S. (1995). The large deviation principle for a general class of queueing systems. I. *Trans. Amer. Math. Soc.* **347** 2689–2751. MR1290716

DUPUIS, P., ELLIS, R. S. and WEISS, A. (1991). Large deviations for Markov processes with discontinuous statistics. I. General upper bounds. *Ann. Probab.* **19** 1280–1297. MR1112416

DUPUIS, P., ISHII, H. and SONER, H. M. (1990). A viscosity solution approach to the asymptotic analysis of queueing systems. *Ann. Probab.* **18** 226–255.

FAYOLLE, G., DE LA FORTELLE, A., LASGOUTTES, J. M., MASSOULIE, L. and ROBERTS, J. (2001). Best-effort networks: Modeling and performance analysis via large networks asymptotics. In *Proc. of IEEE INFOCOM'01*. MR1043946

FAYOLLE, G. and LASGOUTTES, J.-M. (2001). Partage de bande passante dans un ruseau: Approches probabilistes. Technical Report 4202, INRIA.

FAYOLLE, G., MITRANI, I. and IASNOGORODSKI, R. (1980). Sharing a processor among many job classes. *J. Assoc. Comput. Mach.* **27** 519–532.

FREIDLIN, M. I. and WENTZELL, A. D. (1984). *Random Perturbations of Dynamical Systems*. Springer, New York. MR722136

IGNATIOUK-ROBERT, I. (2000). Large deviations of Jackson networks. *Ann. Appl. Probab.* **10** 962–1001.

IGNATYUK, I. A., MALYSHEV, V. A. and SHCHERBAKOV, V. V. (1994). The influence of boundaries in problems on large deviations. *Uspekhi Mat. Nauk.* **49** 43–102. MR1789985

IOFFE, A. D. and TIHOMIROV, V. M. (1979). *Theory of Extremal Problems*. North-Holland, Amsterdam. MR528295

DE LA FORTELLE, A. (2001). Large deviation principle for Markov chains in continuous time. *Problemy Peredachi Informatsii* **36** 120–140.

Shwartz, A. and Weiss, A. (1995). *Large Deviations for Performance Analysis*. Chapman and Hall, London. MR1335456

EDF R&D–1  
av. du Gnral de Gaulle  
92141 Clamart Cedex  
France  
e-mail: franck.delcoigne@edf.fr

INRIA-Domaine de Voluceau  
Rocquencourt  
BP 105  
78153 Le Chesnay Cedex  
France  
e-mail: Arnaud.De_La_Fortelle@inria.fr